\documentclass[12pt,a4paper]{article}
\usepackage{amssymb,amsmath,amsthm,color}
\usepackage{graphicx,mcite}
\usepackage{url}
\usepackage{graphicx}
\usepackage{epsfig}
\newcommand{\Ann}{\mbox{Ann}}
\newcommand{\bt}{\begin{Theorem}}
\newcommand{\bprop}{\begin{Proposition}}
\newcommand{\eprop}{\end{Proposition}}
\newcommand{\et}{\end{Theorem}}
\newcommand{\bl}{\begin{Lemma}}
\newcommand{\el}{\end{Lemma}}
\newcommand{\br}{\begin{Remark}}
\newcommand{\er}{\end{Remark}}

\newcommand{\bi}{\begin{itemize}}
\newcommand{\ei}{\end{itemize}}
\newcommand{\bea}{\begin{eqnarray}}
\newcommand{\eea}{\end{eqnarray}}

\newtheorem{Definition}{Definition}[section]
\newtheorem{Theorem}[Definition]{Theorem}
\newtheorem{Lemma}[Definition]{Lemma}

\newtheorem{Proposition}[Definition]{Proposition}

\newtheorem{Remark}[Definition]{Remark}

\newcommand{\be}{\begin{equation}}
\newcommand{\ee}{\end{equation}}

\theoremstyle{plain} \theoremstyle{plain}
\newtheorem{theorem}{Theorem}[section]
\newtheorem{lemma}[theorem]{Lemma}

\theoremstyle{definition}

\theoremstyle{plain}
\begin{document}
\begin{center}
{\large \bf { Spherical means in annular regions in the $n$-dimensional real hyperbolic spaces}}
\end{center}

\begin{center}
{\bf Rama Rawat  and  R. K. Srivastava}
\end{center}

\begin{abstract}
Let $Z_{r,R}$ be the class of all continuous functions $f$ on the
annulus $\Ann(r,R)$ in the real hyperbolic space $\mathbb B^n$ with
spherical means $M_sf(x)=0$, whenever $s>0$ and $x\in\mathbb B^n$
are such that the sphere $S_s(x)\subset \Ann(r, R) $ and
$B_r(o)\subseteq B_s(x).$ In this article, we give a
characterization for functions in $Z_{r,R}$. In the case $R=\infty$,
this result gives a new proof of Helgason's  support theorem for
spherical means in the real hyperbolic spaces.

\end{abstract}

\section{Introduction}\label{section1}
Let $g$ be a continuous function on the open annulus $\{x\in\mathbb R^d:r<|x|<R\},$ where 
$0\leq r<R\leq\infty$ and $d\geq2.$  We say that $g$ satisfies the Vanishing Spherical Means Condition 
if \[\int_{|x-y|=s}g(y)d\sigma_s(y)=0\] for every sphere
$\{y\in\mathbb R^d: |x-y|=s\}$ which is contained in the annulus and is such that the closed ball 
$\{y\in\mathbb R^d: |y|\leq r\}$ is contained in the closed ball $\{y \in \mathbb R^d: |x-y| \leq s\}.$ 
Here $d\sigma_s$ is the surface measure on the sphere $\{y \in \mathbb R^d: |x-y| =s\}.$

For a continuous function $g$ on $\mathbb R^d,$ let
\begin{equation}\label{exp2}
g(x)=\sum_{k=0}^{\infty}\sum_{j=1}^{d_k}~a_{kj}(\rho)~Y_{kj}(\omega)
\end{equation}
be the spherical harmonic expansion, where $x=\rho\omega,~\rho=|x|,~\omega\in S^{d-1}$
and $\{Y_{kj}(\omega):j=1,\ldots,d_k\}$ is an orthonormal basis for the space $V_k$ of 
homogeneous harmonic polynomials in $d$ variables of degree $k$ restricted to the unit sphere $S^{d-1}.$
Then the following interesting result, which can be thought of as a null space characterization
of the Radon transform over spheres for continuous functions, has been proved in \cite{EK} by 
Epstein and Kleiner:
\begin{theorem}\label{th5} Let $g$ be a continuous function on the annulus
$\{x\in \mathbb R^d :r<|x|<R\},~ 0\leq r<R\leq\infty.$ Then $g$ satisfies the 
Vanishing Spherical Means Condition if and only if
\[a_{kj}(\rho)=\sum_{i=0}^{k-1}~\alpha_{kj}^{i}~\rho^{k-d-2i},\quad \alpha_{kj}^{i}\in\mathbb C,\] 
for all $k\geq1,$ $1\leq j\leq d_k$ and $a_0(\rho)=0,$ whenever $~r<\rho <R.$
\end{theorem}

This result was first proved by Globevnik \cite{G} for the case $n=2.$ For other related work 
we refer to \cite{CQ,Q,V1,V2}.

In a recent work of \cite{NT} the authors have proved a spectral Paley-Wiener theorem for the 
Heisenberg group by means of a support theorem for the twisted spherical means on $\mathbb C^n.$ 
The support theorem for the twisted spherical means can be thought of as a special case of the 
problem described above for the twisted spherical means. The full analogue of the Theorem \ref{th5} 
for these means has been investigated by authors in \cite{RS}.

In this paper, we have investigated the following analogous problem for spherical means in real 
hyperbolic spaces. Let  $\mathbb B^n=\{x\in\mathbb R^n:~ |x|^2=\sum x_i^2<1\}$ be the open unit 
ball in $\mathbb R^n,~n\geq2$, endowed with the Poincare metric $ds^2=\lambda^2(dx_1^2+\cdots+dx_n^2)$, 
where $\lambda=2(1-|x|^2)^{-1}$. Let $B_s(o)=\{x\in\mathbb B^n:d(x,o) \leq s\}$ be the closed geodesic
ball of radius $s$ with centre at the origin and $\Ann(r, R)=\{x\in \mathbb B^n :r<d(x,o)<R\},~0\leq r<R\leq\infty,$
be an open annulus in $\mathbb B^n.$

\bigskip

For $s>0,$ let $\mu_s$ denote the surface measure on the geodesic sphere $S_s(x)=\{y\in\mathbb B^n: d(x,y)=s\}.$
Let $f$ be a continuous function on $\mathbb B^n.$ Define the spherical means of $f$ by
\begin{equation}\label{exp1}
M_sf(x)=\frac{1}{A(s)}\int_{S_s(x)}f(y) d\mu(y),~ x\in \mathbb B^n,
\end{equation}
where $A(s)=(\Omega_n )^{-1}(\sinh s)^{-n+1}.$

\bigskip

Let $Z_{r, R}$ be the class of all continuous functions on $\Ann(r, R)$ with the spherical means 
$M_sf(x)=0,$ whenever $s>0,$ and $x\in\mathbb B^n$ are such that the sphere $S_s(x)\subset\Ann(r, R)$ 
and ball $B_r(o)\subseteq B_s(x).$

\bigskip

Our main result is the following characterization theorem.
\begin{theorem}\label{th1} Let $f$ be a continuous function on $\Ann(r, R)$. 
Then a necessary and sufficient condition for $f$ to be in $Z_{r, R}$ is that 
its spherical harmonic coefficients $a_{kj}(\rho)$ satisfy
\[a_{kj}(\rho)=\sum_{i=1}^{k}C_{kj}^i\frac{(1-\rho^2)^{n+i-2}}{\rho^{n+k-2}},~C_{kj}^i\in\mathbb C,\]
for all $k\geq1,$ $1 \leq j \leq d_k$ and $a_0(\rho)\equiv0,$ whenever  
$~\tanh \frac{r}{2}<\rho<\tanh \frac{R}{2}$.
\end{theorem}

\bigskip

As the authors in \cite{EK} have observed, their result for Euclidean spherical means, 
can be used to derive result for some cases, real hyperbolic spaces being one of them. 
In a recently published book \cite {V2}, Chapter 10, the authors have developed a
general theory of transmutation operators from which our result Theorem \ref{th1} can be 
derived using a string of lemmas proved in this chapter. However, the expressions for the 
spherical harmonic coefficients $a_{kj}(\rho)$  as given in Theorem \ref{th1} have nowhere 
been listed before. Moreover, our approach to the proof of theorem  is direct, transparent 
and brings out the underline geometry of the real hyperbolic spaces clearly. The proof of 
the necessary part of this theorem is close to the work in \cite{V1}, p.108, on problems 
related to spherical means and the proof of the sufficient part is completely new.

The case of other real rank one symmetric spaces can be dealt with in a similar way.

\section{Notation and Preliminaries}\label{section2}

We begin with the realization of real hyperbolic spaces (see \cite{M}, \cite{Re}).
Let $O(1,n+1)$ be the group of all linear transformations which preserve
the quadratic form
\[\langle y,y\rangle=y_0^2-\sum_{i=1}^{n+1}y_i^2, y=(y_0,y_1, \ldots, y_{n+1})\]
on $\mathbb R^{n+2}.$ This group is known as the Lorentz group and is equal to 
\[\{g\in M_{n+2}(\mathbb R):g^tJg=J,~J=\text{diag}(1,-1,\ldots,-1)\}.\]

In particular, $O(1,n+1)$ leaves invariant the cone
 \[C=\left\{y\in\mathbb R^{n+2}: \langle y,y\rangle=0\right\}.\]

With the inhomogeneous coordinates $\eta_i=y_i/ y_0~,i=1,\ldots,n+1,$ the relation 
$\langle y,y\rangle=0$ would imply that  $\eta$ is in $S^n=\{\eta\in\mathbb R^{n+1}:~|\eta|=1\}.$ 
Thus a point on $C$ gets identified with a point on the sphere $S^n.$ Conversely for 
$\eta \in S^n, $ $\eta^*=(1,\eta_1,\ldots,\eta_{n+1})$ gives  a point on the cone $C.$
As $g\in O(1,n+1)$ acts on $\eta^*$ and $g\eta^*\in C,$ $g$ acts on $S^n$ via the above 
identification. More explicitly, $g\eta^*$ can be identified with the point 
$\left(\dfrac{(g\eta^*)_1}{(g\eta^*)_0}, \ldots,  \dfrac{(g\eta^*)_{n+1}}{(g\eta^*)_0}\right)$ in 
$S^n.$ ($(g\eta^*)_0$ is nonzero, as $\eta^*$ is nonzero and $g\eta^*\in C.$ )

\bigskip

Let $O_{\pm}(1,n+1)\cong O(1,n+1)/\{\pm I\}$ be the subgroup of $O(1,n+1)$ which leaves invariant the positive cone
\begin{center}
$C^+=\left\{y=(y_0,y_1,\ldots,y_{n+1})\in\mathbb R^{n+2}: \langle y,y\rangle=y_0^2-\sum_{i=1}^{n+1} y_i^2>0, ~y_0>0 \right\}.$
\end{center}
Equivalently, \[O_{\pm}(1,n+1)=\{g\in M_{n+2}(\mathbb R):g^tJg=J,~J=\text{diag}(1,-1,\ldots,-1),~ g_{00}>0\},\]
where $g_{00}$ is the top left entry in the matrix of $g.$ In particular, $O_{\pm}(1,n+1)$ 
leaves the cone $C^0=\left\{y\in\mathbb R^{n+2}: \langle y,y\rangle=0,~y_0>0\right\}$ invariant.
Moreover, as the action of $g$ and $-g$ in $O(1,n+1)$ on the sphere $S^n$ coincides,
$O_{\pm}(1,n+1)$ also acts on $S^n.$ In fact, this is the group of Mobius transforms
on $S^n$. The real hyperbolic space $\mathbb B^n$ is then isomorphic to the quotient space
~$SO_{\pm}(1,n)/SO(n).$ This isomorphism is established as follows.

\bigskip

We identify $S^n\setminus \{e_{n+1}\}$ with $\mathbb R^n$ under the stereographic projection
from the point $e_{n+1}=(0,\ldots,0,1)\in \mathbb R^{n+1}$ onto the plane $\eta_{n+1}=0.$
Then the $O_{\pm}(1,n+1)$ action on $S^n$ induces an action on $\mathbb R^n\cup \{\infty\}$ 
and vice versa. It turns out that the subgroup of $O_{\pm}(1,n+1)$  which stabilizes
$\mathbb B^n$ is isomorphic to $O_{\pm}(1,n)$. This can be seen as follows.

\bigskip

Let  $x=(x_1,\ldots,x_n)\in \mathbb B^n$. Then the inverse stereographic projection of 
$\eta\in S^n$ of $x$ is given by
\begin{equation}\label{exp9}
\eta_i=\frac{2x_i}{1+|x|^2},~i=1,\ldots,n ~\mbox{and}~\eta_{n+1}=\frac{|x|^2-1}{|x|^2+1}.
\end{equation}
Therefore, $x\in\mathbb B^n $ if and only if $\eta_{n+1}<0.$ Thus a subgroup of $O_{\pm}(1,n+1)$ 
stabilizes the open unit ball $\mathbb B^n$ if and only if it stabilizes the lower hemisphere
$\{\eta\in S^n:~\eta_{n+1}<0\}$. This subgroup in turn is isomorphic to $O_{\pm}(1,n)$, (see \cite{M}).
The elements of this subgroup realized as elements of $O_{\pm}(1,n+1)$ look like
\[\left(
  \begin{array}{cccc}
  g & 0 \\
      0& 1\\
  \end{array}
\right),\]
with $g \in O_{\pm}(1,n).$
Moreover, this action of $O_{\pm}(1,n)$ on $\mathbb B^n$ is transitive and
the orthogonal group $O(n)$ thought of as
 \[\left(
  \begin{array}{cccc}
  1 &  0  \\
  0 &  g  \\
  \end{array}
\right)\]
inside $O_{\pm}(1,n)$ is the isotropy subgroup of the point origin in the ball $\mathbb B^n.$
Thus $\mathbb B^n$ is isomorphic to the quotient space~$O_{\pm}(1,n)/O(n).$ Likewise, 
$\mathbb B^n \cong SO_{\pm}(1,n)/SO(n).$ Let $G=SO_{\pm}(1,n)$ and $K=SO(n)$.Hence onwards, 
we will work with the representation $G/K$ of $\mathbb B^n.$ Using the $G$-invariant metric  
$dy_0^2-dy_1^2-\cdots-dy_n^2$ on the positive cone $y_0^2-y_1^2-\cdots-y_n^2=1,~y_0>0$, $\mathbb B^n$ 
can be endowed with a $G$-invariant Riemannian  metric given by $ds^2= \lambda^2|dx|^2.$
The distance $d(x,y)$ between points $x,y\in \mathbb B^n$, in this metric, is then given by
\[\tanh \frac{1}{2}d(x,y)=\frac{|x-y|}{\sqrt{1-2x.y+|x|^2|y|^2}}.\]
This makes $(\mathbb B^n,d)$ into a Riemannian symmetric space. Group theoretically, $\mathbb B^n=G/K$
is a real rank one symmetric space.

\bigskip

Further, let $G=K\overline{A_+}K$ be the Cartan decomposition of $G$, where
\[A=\small\left\{\left(
    \begin{array}{ccc}
\cosh\frac{t}{2} & 0 & \sinh\frac{t}{2} \\
0 & I_{n-1} & 0 \\
\sinh\frac{t}{2} & 0 &\cosh\frac{t}{2}\\
\end{array}
  \right): ~t\in\mathbb R\right\},
\]
is a maximal abelian subgroup of $G$ and $A_+$ is a chosen positive Weyl chamber $\{a_t : t>0\}.$
Let M be the centralizer $\{k\in K: ka=ak, \forall a\in A\}$ of $A$ in $K$. Therefore, 
$M$ is given by
\[M=\small\left\{\left(
    \begin{array}{ccc}
      1 & 0 & 0 \\
      0 & m &  0 \\
      0& 0 & 1 \\
    \end{array}
  \right)
: ~m\in SO(n-1)\right\}
.\]
Thus the boundary $S^{n-1}$ of $\mathbb B^n $  gets identified with $K/M$ under the map 
$\sigma M\rightarrow \sigma.e_n,$ $\sigma\in K$ where $e_n=(0,0,\ldots,1)\in\mathbb R^n$ 
and the elements of $G/K$ can be thought of as pairs $(a_t,\omega), ~t\geq0,~\omega\in S_{n-1}.$
The point  $(a_t,\omega)$ then is identified with the point $(\cosh\frac{t}{2},~\sinh\frac{t}{2}.\omega)$ 
on the positive cone in $\mathbb R^{n+1}$ and this point in turn, is identified with the point 
$\tanh\frac{t}{2}\omega$ in  $\mathbb B^n$.

\bigskip

Next, we recall certain standard facts about spherical harmonics, for more details see \cite{T}, p. 12.

Let $\hat{K}_M$ denote the set of all the equivalence classes of irreducible unitary representations of 
$K$ which have a nonzero $M$-fixed vector.  It is well known that each representation in $\hat{K}_M$ has 
in fact a unique nonzero $M$-fixed vector, up to a scalar multiple.

For a $\delta \in \hat{K}_M,$ which is realized on $V_{\delta},$ let $\{e_1,\ldots ,e_{d(\delta)}\}$ be an 
orthonormal basis of $V_{\delta},$ with $e_1$ as the $M$-fixed vector. Let 
$t^{ji}_{\delta}(\sigma)=\langle e_i,\delta(\sigma)e_j \rangle ,$ $\sigma\in K$ and $\langle,\rangle $ 
stand for the innerproduct on $V_{\delta}.$ By Peter-Weyl theorem, it follows that
$\{\sqrt{d(\delta)}t^{1j}_{\delta}:1\leq j\leq d(\delta ),\delta\in\hat{K}_M\}$ is
an orthonormal basis of $L^2(K/M).$

We would further need a concrete realization of the representations in $\hat{K}_M,$ which can be done 
in the following way.

Let $\mathbb Z_+$ denote the set of all non-negative integers. For $k\in \mathbb Z_+$, let $P_k$ denote 
the space of all homogeneous polynomials $P$ in $n$ variables of degree $k.$ Let $H_k=\{P\in P_k: \Delta P=0\},$ 
where $\Delta$ is the standard Laplacian on $\mathbb R^n.$ The elements of $H_k$ are called solid 
spherical harmonics of degree $k.$ It is easy to see that the natural action of $K$ leaves the space 
$H_k$ invariant. In fact the corresponding unitary representation $\pi_k$ is in $\hat{K}_M.$ Moreover,
$\hat{K}_M$ can be identified, up to unitary equivalence, with the collection $\{\pi_k: k\in\mathbb Z_+.\}$

\bigskip

Define the spherical harmonics on the sphere $S^{n-1}$ by $Y_{kj}(\omega )=\sqrt{d_k}t^{1j}_{\pi_k}(\sigma),$ 
where $\omega =\sigma.e_n\in S^{n-1},$ $\sigma\in K$ and $d_k$ is the dimension of $H_k.$  Then 
$\{Y_{kj}:1\leq j\leq d_k,k\in\mathbb Z_+\}$ forms an orthonormal basis for $L^2(S^{n-1}).$
Therefore, for a continuous function $f$ on $\mathbb B^n,$ writing
$y=\rho \,\omega,$ where $0<\rho<1$ and $\omega \in S^{n-1},$ we can expand the function $f$ in terms of 
spherical harmonics as in the (\ref{exp2}). For each non negative integer $k$, the $k^{th}$ spherical 
harmonic projection, $\Pi_k(f)$ of the function $f$ is defined by \begin{equation}\label{exp11}
\Pi_k(f)(y)= ~\sum_{j=1}^{d_k}~a_{kj}(\rho)~Y_{kj}(\omega),
\end{equation}
where $a_{kj}$'s are the spherical harmonic coefficients of the function $f.$

\section{Auxiliary results}\label{section3}

We begin with the observation that the $K$-invariance of the annulus and the measure $\mu_s$ implies that 
for any $f$ in $Z_{r,R}$ and $k \in \mathbb Z_+$, $\Pi_k(f),$ as defined in equation (\ref{exp11}), also 
belongs to $Z_{r,R}.$ In fact the following stronger result is true.

\begin{lemma}\label{lemma1} Let $f \in Z_{r, R}$. Then each of the spherical harmonic projection 
$\Pi_k(f)\in Z_{r,R},$ in fact $a_{kj}(\rho)Y_{km}(\omega)\in Z_{r, R}~\forall~j,m,
~1\leq j,m\leq d_k$ and $\forall ~k\geq0.$
\end{lemma}

\begin{proof}
Since the measure $\mu_s$ and space $\Ann(r,R)$ both are rotation invariant, 
it is easy to verify that, if $f\in Z_{r, R}$, then the function $f(\tau.y)\in Z_{r, R}$ 
for each $\tau\in K$. Since the space $H_k$ is $K$-invariant, for $\tau\in K$ and 
a spherical harmonic $Y_{kj}$, we have
\[Y_{kj}(\tau^{-1}\omega)=\sum_{m=1}^{d_k}\overline{t^{mj}_{\pi_k}(\tau)}Y_{km}(\omega).\]
Hence from the equation (\ref{exp2}), the function $f(\tau^{-1}.)$ can be decomposed as
\[f(\tau^{-1}\rho \omega)=\sum_{k\geq0}\sum_{j,m=1}^{d_k}a_{kj}(\rho)\overline{t^{mj}_{\pi_k}(\tau)}Y_{km}(\omega).\]
The set $\{\sqrt{d_k}t^{mj}_{\pi_k}: 1\leq j,m \leq d_k, k\geq0\}$ form an orthonomal basis for $L^2(K),$
therefore,
\[a_{kj}(\rho)Y_{km}(\omega)=d_k\int_K f(\tau^{-1}\rho\omega)t^{mj}_{\pi_k}(\tau)d\tau\in Z_{r,R}.\]
Subsequently, each projection $\Pi_k(f)$ belongs to $Z_{r,R}$.
\end{proof}

\bigskip

Next, we need the following explicit expression for the action of $G$ on $\mathbb B^n,$ which has been 
derived in \cite{J}.

\begin{lemma}\label{lemma4}
Let $g\in G$ and  $\;x\in \mathbb B^n$. Then  $g.(x_1,\ldots, x_n)=(y_1,\ldots,y_n)$, where
\begin{equation}\label{exp8}
y_j=\frac{\frac{(1+|x|^2)}{2}g_{j0}+\sum_{l=1}^n g_{jl}x_l}{\frac{1-|x|^2}{2}+\frac{(1+|x|^2)}{2}g_{00}+
\sum_{l=1}^n g_{0l}x_l}, ~j=1,\ldots,n.
\end{equation}
\end{lemma}
\begin{proof}
By equation (\ref{exp9}), a point $x\in\mathbb B^n$ is mapped to the point $\eta\in S^n$ via the
the inverse stereographic projection. By definition, for $g\in G$,
\[g\cdot\eta=\left(
  \begin{array}{cccc}
 g &  0 \\
 0 &   1\\
  \end{array}
\right)\left(
\begin{array}{c}
           1 \\
         \eta \\
\end{array}
       \right)=\alpha
       ,\]
where  $\alpha=(\alpha_0,\ldots,\alpha_n,\eta_{n+1})$ and $\alpha_j= g_{j0}+\sum_{l=1}^ng_{jl}\eta_l,~l=0,1,\ldots,n$.
Since the cone $C^0$ is $G$-invariant, it follows that $\alpha_0>0$. In the inhomogeneous coordinates, 
introduced earlier, the point $\alpha$ gets identified with the point
$\left(\dfrac{\alpha_1}{\alpha_0},\ldots,\dfrac{\alpha_n}{\alpha_0},\dfrac{\eta_{n+1}}{\alpha_0}\right)$
on the sphere $S^n$. The image of this point, under the stereographic projection is the point
$y=(y_1,\ldots, y_n)\in \mathbb B^n$, where
\[y_j=\frac{\alpha_j/\alpha_0}{1-\eta_{n+1}/\alpha_0},~j=1,\ldots,n.\]
That is,
\[y_j=\frac{g_{j0}+\sum_{l=1}^ng_{jl}\eta_l}{g_{00}+\sum_{l=1}^ng_{0l}\eta_l-\eta_{n+1}},~j=1,\ldots,n.\]
Since we know that
\[\eta_l=\frac{2x_l}{1+|x|^2},~l=1,\ldots,n ,~\eta_{n+1}=\frac{|x|^2-1}{|x|^2+1},\]
a simple computation gives
\[y_j=\frac{\frac{(1+|x|^2)}{2}g_{j0}+\sum_{l=1}^n g_{jl}x_l}{\frac{1-|x|^2}{2}+\frac{(1+|x|^2)}{2}g_{00}+
\sum_{l=1}^n g_{0l}x_l}, ~~j=1,\ldots,n.\]
\end{proof}

\bigskip

As in the proof of the Euclidean case \cite{EK}, to characterize functions in $Z_{r, R}$ it would be enough 
to characterize the spherical harmonic coefficients of smooth functions in $Z_{r,R}$. This can be done using 
the following approximation argument.

Let $\varphi_\epsilon$ be nonnegative, $K$-biinvariant, smooth, compactly supported approximate identity 
on $G/K$. Let $f\in Z_{r,R}$. Then $f$ can be thought of as a right $K$-invariant function on $G$. Define
\[S_\epsilon(f)(g)=\int_G f(gh^{-1})\varphi_\epsilon(h)dh,~ g\in G.\] Then $S_\epsilon(f)$ is smooth and 
it is easy to see that $S_\epsilon(f)\in Z_{r+\epsilon, R-\epsilon}$ for each $\epsilon >0.$ Since $f$ is
continuous, $S_\epsilon(f)$ converges to $f$ uniformly on compact sets. Therefore, for each $k$,
\[\lim\limits_{\epsilon \rightarrow 0}\Pi_k(S_\epsilon(f))=\Pi_k(f).\]
Hence, we can assume, without loss of generality, that the functions in $Z_{r,R}$ are also smooth in the 
annulus $\Ann(r, R).$

\bigskip

We next introduce right $K$-invariant differential operators on $G$ which leave invariant the space 
$Z_{r,R}$. These differential operators arise naturally from the Lie algebra $\mathfrak{g}$ of $G$, 
in the following way. They also appear prominently in the work of  Volchkov on ball means in real
hyperbolic spaces, (\cite{V1}, p. 108).

\bigskip

Let $\mathfrak{g}=\mathfrak{k}+\mathfrak{p}$ be the Cartan decomposition of the Lie algebra $\mathfrak{g}$ 
of $G$. Here  $\mathfrak{k}$ is the Lie algebra of $K$ and $\mathfrak{p}$ its orthogonal complement in 
$\mathfrak{g}$ with respect to the killing form $B(-,-)$. Let $X_i=E_{0i}+ E_{i0}, ~i=1,\ldots,n$ and
$X_{ij}=E_{ij}-E_{ji}, ~1\leq i<j\leq n$, where $E_{ij}\in gl(n+1,\mathbb R)$ is the matrix with entry 
$1$ at the ${ij}^{th}$ place and zero elsewhere. Then  $\{X_i:i=1,\ldots,n\}$ and $\{X_{ij}:1\leq i<j\leq n\}$ 
form bases of $\mathfrak{p}$ and $\mathfrak{k}$ respectively.

\bigskip

Let $f\in C^{\infty}(\mathbb B^n)$. Then $f$ can be thought of as the right $K$-invariant function on $G$.
For a given $X\in \mathfrak{g}$, let $\tilde X$ be the differential operator given by
\begin{eqnarray}\label{exp7}
(\tilde{X}f)(gK)=\left.\frac{d}{dt}\right\rvert_{t=0} f(\exp tXgK).
\end{eqnarray}
For $X=X_p\in \mathfrak{p}$, let
\[\tau_{t,p}=\exp tX_p=\small\left(\begin{array}{ccccccccc}
\cosh t & 0&\sinh t& 0 \\
0 & I_{p-1} &0  & 0 \\
\sinh t& 0 & \cosh t & 0 \\
0 & 0 & 0 & I_{n-p}
\end{array}\right),\]
for $t \in \mathbb R$. Let $x\in\mathbb B^n$. Then by Lemma \ref{lemma4},
$\tau_{t,p}.x=y\in\mathbb B^n,$ where
$y_j=x_ju(t,x), ~\mbox{if}~ j\neq p $ and $y_p=(x_p\cosh t+(1+|x|^2)\frac{\sinh t}{2}) u(t,x)$,
 $ u(t,x)=(\cosh^2\frac{t}{2}+x_p\sinh t+ |x|^2\sinh^2\frac{t}{2})^{-1}$.
Rewrite $\tau _{t,p}.x$ as $\tau(t, x)$. Then $\tau$ is a differentiable function on $\mathbb R\times\mathbb R^n$
 into $\mathbb R^n$ and from
(\ref{exp7}), we have
\[\frac{\partial}{\partial t}(fo\tau(t,x))=f'(\tau(t,x))\frac{\partial\tau}{\partial t}(t,x)
=\sum_{j=1}^n\frac{\partial f}{\partial y_j}\frac{\partial y_j}{\partial t}.\]
Evaluating the above equation at $t=0$, we get
\begin{eqnarray}\label{exp5}
\left.\frac{\partial}{\partial t}(fo\tau(t,x))\right\rvert_{t=0}=\sum_{j=1}^n\left.\frac{\partial f}{\partial y_j}\right\rvert_{t=0} 
\left.\frac{\partial y_j}{\partial t}\right\rvert_{t=0}=
\sum_{j=1}^n\frac{\partial f}{\partial x_j}\left.\frac{\partial y_j}{\partial t}\right\rvert_{t=0}.
\end{eqnarray}
A straightforward calculation then gives,
\[
\left.\frac{\partial y_j}{\partial t}\right\rvert_{t=0}=\left\{
\begin{array}{ll}
- x_px_j, ~~~~~~~~~~~~~~~~~~~\mbox{if}~ j\neq p;\\
\frac{1}{2}(1+|x|^2)-x_p^2,~~~~~~~\textrm{if}~j=p.
\end{array}\right.\]
Substituting these values in (\ref{exp5}), we get
\[\tilde{X}_p =\frac{1}{2}(1+|x|^2)\frac{\partial}{\partial x_p}-\sum_{j=1}^n~x_px_j\frac{\partial}{\partial x_j},~p=1,\ldots,n.\]

The following lemma is a crucial step towards the proof of our main result.
\begin{lemma}\label{lemma2}
Suppose $f$ is a smooth function belonging to $Z_{r,R}$. Then the function $\tilde{X}_pf\in Z_{r,R}, ~\forall~p,~1\leq p\leq n.$
 \end{lemma}

\begin{proof}
For a fixed $t\in\mathbb R$, define
\[\epsilon_{1}=\sup_{y\in B_r(o)}d(\tau_{t,p}.y,y)\text{ and }\epsilon_{2}=\sup_{y\in B_R(0)}d(\tau_{t,p}.y,y).\]
Then, it is easy to see that the translated function $\tau_{t,p}f$ defined by 
$\tau_{t,p}f(y)=f(\tau_{t,p}.y),~y\in\mathbb B^n$ belongs to 
$Z_{r+\epsilon_{1}, R-\epsilon_{2}}.$ Therefore,
\[\int_{S_{s}(x) }~f(\tau _{t,p}.\xi)d\mu_{s}(\xi)=\int_{S_{s}(\tau _{t,p}.x) }~f(\xi)d\mu_{s}(\xi)=0,\]
whenever $S_s(x)\subset\Ann(r+\epsilon_{1}, R-\epsilon_{2})$ and
$B_{r+\epsilon_{1}}(0)\subset B_s(x).$ As $t\rightarrow0$, this implies
\[\int_{S_{s}(x) }\left.\frac{\partial f}{\partial t}\right\rvert_{t=0}(\tau _{t,p}.\xi)d\mu_{s}(\xi)=0,\]
whenever $S_s(x)\subset\Ann(r, R)$ and $B_{r}(0)\subseteq B_s(x).$
Hence $\tilde{X}_pf\in Z_{r,R}.$
\end{proof}

A repeated application of Lemma \ref{lemma2}, leads naturally to a family of differential 
operators which we now introduce. These operators also appear in the work of Volchkov 
(\cite{V1}, p.108) in the problems on averages over geodesic balls in real hyperbolic
spaces. Let $C^1(0,1)$ denote the space of all differentiable functions on $(0,1)$. For 
$m \in \mathbb Z$, the set of integers, define a differential operator $A_m$ on  $C^1(0,1)$ by
\begin{equation}\label{exp25}
 \mathbb (A_mf)(t):=\frac{t^m}{(1-t^2)^{m-1}} \frac{d}{dt} \left[\left(\frac{1}{t}-t \right)^m f(t) \right].
\end{equation}
The Laplace-Beltrami operator $\mathcal L_x$ on $\mathbb B^n$ (\cite{H}, p.31)
 is given by
\[\mathcal L_x=\frac{(1-|x|^2)^n}{4}\sum_i\frac{\partial}{\partial x_i}\left(\sum_i(1-|x|^2)^{2-n}\frac{\partial}{\partial x_i} \right).\]
The radial part $\mathcal L_s$ of $\mathcal L_x$ is given by
\[\mathcal L_s=\frac{\partial^2}{\partial s^2}+(n-1)\coth s \;\frac{\partial}{\partial s}\]
and satisfies the Darboux equation $M_s\mathcal L_x=\mathcal L_sM_s$ (\cite{H}, p.159). For any positive integer $k$, 
define an operator $\mathcal L_k$ by
\begin{equation}\label{exp3}
\mathcal L_k=\mathcal L_x-4(k-1)(n+k-2)\text{Id}.
\end{equation}
Let $f(x)=a(\rho)Y_k(\omega),$ where $Y_k$ is a spherical harmonic of degree $k$.
Then, a simple calculation shows that
\[\mathcal L_kf(x) =(A_{k-1}A_{2-k-n}a)(\rho)Y_k(\omega),~x=\rho\omega.\]
\begin{lemma}\label{lemma3}
Let $x=\rho\omega, ~0<\rho<1 $ and $\omega\in S^{n-1}$ and $k\geq0$. Suppose the function 
$f(x)=a(\rho)Y_k(\omega) \in Z_{r,R}$. Then
\begin{description}
\item[(i)] $(A_{2-k-n}a)(\rho)Y_{(k-1)j}(\omega) \in Z_{r,R}, k\geq1$ and $1\leq j\leq d_{k-1}(n),$

\item[(ii)] $(A_k a)(\rho)Y_{(k+1)i}(\omega)\in Z_{r,R}, k\geq0$ and $1\leq i\leq d_{k+1}(n),$

\item[(iii)] $(A_{1-k-n}A_k a)(\rho)Y_k(\omega)\in Z_{r,R}, k\geq0$ and

\item[(iv)] $\mathcal L_kf(x) =(A_{k-1}A_{2-k-n}a)(\rho)Y_k(\omega)\in Z_{r,R}, ~k\geq1.$
\end{description}
\end{lemma}
\begin{proof}
Let $k\geq1$. Let $P(x)=\rho^k Y_k(\omega)$ and $\tilde{a}(\rho)=\rho^{-k}a(\rho)$. Then $f=\tilde{a}P$,
where $P\in H_k$. By Lemma \ref{lemma2}, the function $2\tilde{X_p}f\in Z_{r,R}$ $\forall ~1\leq p\leq n.$ 
A straightforward calculation then gives
\begin{eqnarray}\label{exp6}
2\tilde{X_p}f=\left(\frac{(1-\rho^2)}{\rho}\frac{\partial\tilde{a}}{\partial\rho}-2k\tilde{a}\right)x_pP +(1+\rho^2)
\tilde a\frac{\partial P}{\partial x_p}.
\end{eqnarray}
Further,
\[x_pP=P_{k+1}+\frac{|x|^2}{n+2(k-1)}\frac{\partial P}{\partial x_p},\]
where $P_{k+1} \in H_{k+1}$ (for a proof, see \cite{EK}). Let $l=2-k-n$, then (\ref{exp6}) gives
\[2\tilde{X_p}f=\left(\frac{(1-\rho^2)}{\rho}\frac{\partial\tilde{a}}{\partial\rho}-2k\tilde{a}\right)
\left( P_{k+1}+\frac{\rho^2}{k-l}\frac{\partial P}{\partial x_p}\right) +(1+\rho^2)\tilde a\frac{\partial P}{\partial x_p}.\]
After a rearrangement of terms, we get
\begin{eqnarray*}
2(k-l)\tilde{X_p}f&=&(k-l)\left(\frac{(1-\rho^2)}{\rho}\frac{\partial\tilde{a}}{\partial\rho}-2k\tilde{a}\right)P_{k+1}\nonumber\\
                 &+&\left(\rho(1-\rho^2)\frac{\partial\tilde{a}}{\partial\rho}-2k\rho^2\tilde{a}+(k-l)(1+\rho^2)\tilde a \right)
\frac{\partial P}{\partial x_p}.
\end{eqnarray*}
Since $\tilde a(\rho)=\rho^{-k}~a(\rho)$, $\dfrac{\partial\tilde a}{\partial\rho}=-k\rho^{-k-1}a+
\rho^{-k}\dfrac{\partial a}{\partial\rho}.$ Using this in the above equation, we have
\begin{eqnarray}\label{exp24}
2(k-l)\tilde{X_p}f&=&(k-l)\left((1-\rho^2)\frac{\partial a}{\partial\rho}-k\frac{(1+\rho^2)}{\rho}a\right)\rho^{-k-1}P_{k+1}\nonumber\\
&+&\left((1-\rho^2)\frac{\partial a}{\partial\rho}-l\frac{(1+\rho^2)}{\rho} a\right)\rho^{-k}\frac{\partial P}{\partial x_p}.
\end{eqnarray}
Also the operator $A_m$, given by (\ref{exp25}), can be rewritten as
\[A_m=(1-t^2)\frac{d}{dt}-m\frac{(1+t^2)}{t}.\]
Thus (\ref{exp24}) can be rephrased as
\[2(k-l)\tilde{X_p}f=(A_k a)(\rho)\rho^{-k-1}P_{k+1}+(A_{2-k-n}a)(\rho)\rho^{-k+1}\frac{\partial P}{\partial x_p}\in Z_{r,R},\]
whenever $1\leq p\leq n$. Consequently, by Lemma \ref{lemma1}, we get $(A_k a)(\rho)\rho^{-k-1}P_{k+1}\in Z_{r,R}$ and $(A_{2-k-n}
a)(\rho)\rho^{-k+1}\dfrac{\partial P}{\partial x_p}\in Z_{r,R}$ and in particular,  $(A_{2-k-n}a)(\rho)~Y_{(k-1)j}(\omega)$ 
and $(A_ka)(\rho)Y_{(k+1)i}(\omega)$ belong to $Z_{r,R}.$

The assertions (iii) and (iv)  can be obtained by composing (i) and (ii).
\end{proof}

\section{Proof of the main result}\label{section4}
In this section, we prove our main result Theorem \ref{th1}. We first take up the necessary part
of the theorem.
\begin{Proposition}\label{pro1}
Let  $f$ be a radial function in $Z_{r,R}$. Then $f\equiv0$ on $\Ann(r, R)$.
\end{Proposition}
\begin{proof}
By hypothesis \[\int_{S_s(x) }~f(\rho) d\mu_s(y)=0,\] whenever $x\in \mathbb B^n$ is such that the sphere
$S_s(x)\subseteq \Ann(r, R)$ and ball $B_r(o)\subseteq B_s(x).$ Evaluating at $x=0$, this implies
\[\int_{S_s(0)}~f(|y|)d\mu_s(y)=0, \mbox{ whenever } ~~R>s>r.\] Thus $f(\tanh \frac{s}{2})=0, ~R>s>r$.
\end{proof}
\begin{Proposition}\label{pro2} Let $f(\rho \omega)=a(\rho)Y_k(\omega) \in Z_{r,R},k\geq1.$
Then  $a(\rho)$ is given by is given by
\begin{equation}\label{exp26}
a(\rho)=\sum_{i=1}^{k}C_i\frac{(1-\rho^2)^{n+i-2}}{\rho^{n+k-2}}, ~C_i\in \mathbb C,
\mbox{ whenever }~\tanh \frac{r}{2}<\rho<\tanh \frac{R}{2}.
\end{equation}
\end{Proposition}
\begin{proof}
We use induction on $k$. For $k=1$, let $f(\rho\omega)=a(\rho)Y_1(\omega) \in Z_{r,R}$. Using Lemma
\ref{lemma3}({ii}), it follows that  $(A_{1-n}a)(\rho)Y_0(\omega)$ belongs to $ Z_{r,R}$. Therefore, 
by Proposition \ref{pro1}, $(A_{1-n}a)(\rho)=0,$ on $\Ann(r, R).$ On solving this differential
equation, we get $a(\rho)=C\left(\frac{1}{\rho}-\rho \right)^{n-1}.$

\bigskip

Next, we assume the result is true for $k$. Suppose $f(\rho\omega)=a(\rho)Y_{k+1}(\omega) \in Z_{r,R}$. 
An application of Lemma \ref{lemma3}({ii}) gives $A_{1-k-n}a(\rho)Y_k(\omega)\in Z_{r,R}$. Using the 
result for $k$ and the definition of $A_{1-k-n}$, it follows that
\[\frac{\rho^{1-k-n}}{(1-\rho^2)^{-k-n}}\frac{\partial}{\partial\rho}\left(\left(\frac{1}{\rho}-\rho\right)^{1-k-n} a(\rho)\right)=
\sum_{i=1}^{k}C_i\frac{(1-\rho^2)^{n+i-2}}{\rho^{n+k-2}}.\]
Simplifying this equation and integrating both the sides with respect to $\rho$, we obtain
\[\left(\frac{1}{\rho}-\rho\right)^{1-k-n} a(\rho)=\sum_{i=1}^{k}D_i\frac{1}{(1-\rho^2)^{k-i+2}}+D_{k+1}, D_i\in\mathbb C.\] 
Hence
\[a(\rho)=\sum_{i=1}^{k+1}D_i\frac{(1-\rho^2)^{n+i-2}}{\rho^{n+k-1}},\] whenever $~~\tanh\frac{r}{2}<\rho<\tanh\frac{R}{2}$.
\end{proof}
Now, we shall prove the sufficient part of Theorem \ref{th1}. For this, 
without loss of generality, we may assume that $R=\infty$. The idea of 
the proof is to use the asymptotic behavior of the hypergeometric function 
and compare it with that of the coefficients given in (\ref{exp26}). In the 
proof, we need, the following result from Erdelyi et al. \cite{EMOT},~p. 75-76.
\begin{lemma}\label{lemma5}
The general solution of the hypergeometric differential equation
\begin{equation}\label{exp15}
z(1-z)U''+\{\gamma-(\alpha+\beta+1)z\}U'-\alpha\beta U=0,
\end{equation}
where $\alpha, \beta, \gamma$ are independent of $z,$ in the neighborhood of $\infty$ is given in the
following way. If ~$\alpha-\beta$~ is not an integer then
\[U(z)=\lambda_1 z^{-\alpha}+\lambda_2 z^{-\beta}+ O\left(z^{-\alpha-1}\right)+O\left(z^{-\beta-1}\right),\]
where $\lambda_1(x)$ and $\lambda_2(x)$ are non zero. Otherwise,  $z^{-\alpha}$ or $z^{-\beta}$ has to be 
multiplied by a factor of  ~$\log z$.
\end{lemma}

\begin{theorem}\label{th5}
Let $y=\rho\omega, ~\omega\in S^{n-1}$ and $\tanh\frac{r}{2}< \rho<\infty.$ Let $h(y)=a(\rho)Y_k(\omega)$
with
\[a(\rho)=\sum_{i=1}^{k}C_i\frac{(1-\rho^2)^{n+i-2}}{\rho^{n+k-2}},~C_i\in \mathbb C.\]
Then the function $h\in Z_{r, \infty}.$
\end{theorem}

\begin{proof}
We use the induction hypothesis on $k$. For $k=1,$ consider the function
$h(y)=\left(\frac{1}{\rho}-\rho \right)^{n-1}Y_1(\omega).$ In this case, a straightforward
calculation gives $(A_0A_{1-n})\left(\frac{1}{\rho}-\rho\right)^{n-1}\equiv0.$ 
Thus, we can write
\[\mathcal L_1h(y)=(A_0A_{1-n})\left(\frac{1}{\rho}-\rho \right)^{n-1}Y_1(\omega)=(A_0A_{1-n})h(y)=0.\]
Hence, from Equation (\ref{exp3}), it follows that $\mathcal L_yh(y)=0.$ Again by Darboux's equation 
$\mathcal L_sM_sh=M_s\mathcal L_yh,$ the above leads to $\mathcal L_s(M_sh)=0.$ Define $F_1(s,x)=M_sh(x)$. 
For fixed $x, F_1$ as a function of $s$ satisfies the differential equation
\begin{equation}\label{exp12}
\frac{\partial^2F_1}{\partial s^2}+(n-1)\coth s~\frac{\partial F_1}{\partial s}=0.
\end{equation}
Setting $z=-\sinh^2\frac{s}{2},$ we get
\[ \frac{\partial F_1}{\partial s}=\frac{\partial F_1}{\partial z}\frac{\partial z}{\partial s}
= -\frac{1}{2}\sinh s\frac{\partial F_1}{\partial z},~
 \frac{\partial^2 F_1}{\partial s^2}=\frac{1}{4}{(\sinh s)}^2\frac{\partial^2 F_1}{\partial z^2}
-\frac{1}{2}\cosh s\frac{\partial F_1}{\partial z}.\]
After substituting these values in (\ref{exp12}), we obtain
\begin{equation}\label{exp13}
-z(1-z)\frac{\partial^2 F_1}{\partial z^2}-\left(\frac{n}{2}-nz\right)\frac{\partial F_1}{\partial z}=0.
\end{equation}
Comparing this equation with (\ref{exp15}), we get $\gamma=\frac{n}{2}, \alpha+\beta+1=n, \alpha\beta=0$.
For $\alpha=0$, $\beta=n-1$. The solution of (\ref{exp13}) as $|z|\rightarrow\infty$ is given by
\begin{equation}\label{exp16}
F_1(z,x)=\lambda_1(x)+\lambda_2(x) z^{-(n-1)}\log z+ O\left(z^{-(n-1)-1}\log z\right),
\end{equation}
where $\lambda_1(x)$ and $\lambda_2(x)$ are non zero. On the other hand, for $x=g.o, g\in G,$ we have
\begin{eqnarray*}
M_sh(x)&=&\frac{1}{A(s)}\int_{S_s(x)}h(y) d\mu_s(y),\\
      &=&\frac{1}{A(s)}\int_{S_s(o)}h(g.y) d\mu_s(y).
\end{eqnarray*}
From the above equation, it follows that
\begin{equation}\label{exp18}
M_sh(x)=O\left(a\left(\tanh\frac{s}{2}\right)\right), \mbox{ as } s\rightarrow\infty.
\end{equation}
From (\ref{exp18}) one can also conclude that any function of type $h(y)=a(\rho)Y_k(\omega),$ 
must satisfy the relation $M_sh(x)=O\left(a\left(\tanh\frac{s}{2}\right)\right)$.
In fact, for $k=1$,
\begin{equation}\label{exp22}
\left|a\left(\tanh\frac{s}{2}\right)\right|=\left|\cosh\frac{s}{2}\sinh\frac{s}{2}\right|^{-(n-1)}
=\left|z(1-z)\right|^{-\frac{(n-1)}{2}}.
\end{equation}
From (\ref{exp18}) and (\ref{exp22}), we have $F_1(z,x)=O(z^{-(n-1)})$, as $|z|\rightarrow\infty$.
In view of (\ref{exp16}), we infer that $F_1(z,x)=0$, whenever $|z|>\sinh^2r$. Thus $M_sh(x)=0$, 
whenever $x\in \mathbb B^n$ is such that the ball $B_r(o)\subseteq B_s(x)$ and $r<s<\infty$, 
which proves the result for $k=1$.

To complete the induction argument, we assume the result is true for $k-1$ and
then prove for $k$. For this, consider the function
\[h(y)=a(\rho)Y_k(\omega)=\frac{(1-\rho^2)^{n+i-2}}{\rho^{n+k-2}}Y_k(\omega),\]
for each $i,~1\leq i\leq k.$ Using Lemma \ref{lemma3}(i) and the case $(k-1)$, it follows that
\[(A_{2-k-n}a)(\rho)Y_{k-1}(\omega)=\frac{(1-\rho^2)^{n+i-2}}{\rho^{n+k-3}}Y_{k-1}(\omega)\in Z_{r,\infty}.\]
Applying Lemma \ref{lemma3}(ii), it follows that $\mathcal L_kh(y)=(A_{k-1}A_{2-k-n})a(\rho)Y_k(\omega)$ belongs 
to $Z_{r,\infty}.$ Since we know that
\[\mathcal L_kh(y)=\mathcal L_yh(y)-4(k-1)(n+k-2)h(y),\] 
therefore, evaluating mean and using Darboux's equation, we obtain
\[\mathcal L_s(M_sh(x))-4(k-1)(n+k-2)M_sh(x)=0,\]  
whenever $x\in \mathbb B^n$ is such that the ball $B_r(o)\subseteq B_s(x)$ and $r<s<\infty$. 
Let $F_k(s,x)=M_sh(x)$. For fixed $x, F_k$ as a function of $s$ satisfies the differential equation
\[\frac{\partial^2F_k}{\partial s^2}+(n-1)\coth s~\frac{\partial F_k}{\partial s}-4(k-1)(n+k-2)F_k=0.\]
Using the change of variable $z=-\sinh^2\frac{s}{2}$, the above equation becomes
\begin{equation}\label{exp29}
-z(1-z)\frac{\partial^2 F_k}{\partial z^2}-\left(\frac{n}{2}-nz\right)\frac{\partial F_k}{\partial z}-4(k-1)(n+k-2)F_k=0.
\end{equation}
Comparing this equation with (\ref{exp15}), we have $\gamma=\frac{n}{2},\alpha+\beta+1=n, \alpha\beta=-4(k-1)(n+k-2)$. 
On solving, we find that $\alpha=\frac{n-1+\nu}{2},~\beta=\frac{n-1-\nu}{2},$ where $\nu=\sqrt{{(n-1)}^2+16(k-1)(n+k-2)}.$ 
If $\nu\not\in\mathbb Z,$ then the solution of (\ref{exp29}) as $|z|\rightarrow\infty$
is given by
\begin{equation}\label{exp30}
F_k(z,x)=\lambda_1(x) z^{-\alpha}+\lambda_2(x) z^{-\beta}+ O\left(z^{-\alpha-1}\right)+O\left(z^{-\beta-1}\right),
\end{equation}
where $\lambda_1(x)$ and $\lambda_2(x)$ are non zero, otherwise 
$z^{-\alpha}$ or $z^{-\beta}$ has to be multiplied by a factor of $\log z$.
From the given expression of the function $h$, we find that
\[M_sh(x)=O\left(a\left(\tanh\frac{s}{2}\right)\right),\mbox{ as } s\rightarrow\infty.\]
Using $z=-\sinh^2\frac{s}{2}$, it follows that
\[\left|a\left(\tanh\frac{s}{2}\right)\right|=\left|\frac{\left(\text{sech}^2\frac{s}{2}\right)^{n+i-2}}{\left(\tanh\frac{s}{2}\right)^{n+k-2}}\right|
=\frac{|1-z|^{\frac{n+k-2}{2}-(n+i-2)}}{|z|^{\frac{n+k-2}{2}}}.\]
That is,
\begin{equation}\label{exp4}
F_k(z,x)=O\left(z^{-(n+i-2)}\right), ~i=1,\ldots,k \mbox{ as }|z|\rightarrow\infty,
\end{equation}
which contradicts with the expression of $F_k(z,x)$ given by (\ref{exp30}).
Therefore, $F_k(z,x)=0$, whenever $|z|>\sinh^2r$. Hence, we conclude that $M_sh(x)=0$, 
whenever $x\in\mathbb B^n$ is such that the ball $B_r(o)\subseteq B_s(x)$ and 
$r<s<\infty$, which proves the result for any positive integer $k$. 
This completes the proof.
\end{proof}

As a corollary of Theorem \ref{th1}, we have the following Helgason support theorem (see \cite{H}, p. 156).
\begin{theorem} \label{th4}
Let $f$ be a function on $\mathbb B^n$. Suppose for each  $m \in\mathbb Z_+$, the function 
$e^{md(x,~0)}f(x)$ is  bounded. Then $f$ is supported in closed geodesic ball $B_r(o)$ 
if and only if $f \in Z_{r,\infty}$.
\end{theorem}

\begin{proof} The decay condition on function $f$ implies that for all $k$ and $j$, $a_{kj}(|x|)=0,$
whenever $|x|>\tanh\frac{r}{2}$. This proves $f$ is supported in the ball $B_r(o)$.
\end{proof}

\noindent{\bf Acknowledgements:} We thank the referee for the valuable remarks. 
The second author wishes to thank the MHRD, India, for the senior research 
fellowship and IIT Kanpur for the support provided during the preparation of this work.

\vskip.15in
\begin{flushleft}
Department of Mathematics and Statistics, \\
Indian Institute of Technology \\
Kanpur 208 016, India.\\
E-mail:~rrawat@iitk.ac.in, rksri76@gmail.com\\
\end{flushleft}
\end{document}